\documentclass{amsart}
\usepackage{amsmath}
\usepackage{amstext}
\usepackage{amssymb}
\usepackage{amsthm}
\swapnumbers
\theoremstyle{plain}%default
\newtheorem{thm}{Theorem}[section]
\newtheorem{lem}[thm]{Lemma}
\newtheorem{lemdef}[thm]{Lemma and Definition}
\newtheorem{prop}[thm]{Proposition}
\newtheorem{cor}[thm]{Corollary}

\theoremstyle{definition}
\newtheorem{defi}[thm]{Definition}

\newtheorem{ntn}[thm]{Notation}
\newtheorem{rmd}[thm]{Reminder}

\newtheorem{disc}[thm]{Discussion}
\theoremstyle{remark}

\newtheorem{rmk}[thm]{Remark}

 \DeclareMathOperator{\height}{ht}
 
\DeclareMathOperator{\Ext}{Ext} \DeclareMathOperator{\Hom}{Hom}
\DeclareMathOperator{\Endo}{End} 
 
 \DeclareMathOperator{\depth}{depth}

\DeclareMathOperator{\Ass}{Ass} \DeclareMathOperator{\ass}{ass}

\def\Z{\mathbb Z}
\def\N{\mathbb N}

\def\fa{{\mathfrak{a}}}

\def\fb{{\mathfrak{b}}}

\def\fc{{\mathfrak{c}}}

\def\fm{{\mathfrak{m}}}

\def\fp{{\mathfrak{p}}}

\def\fq{{\mathfrak{q}}}

\def\fs{{\mathfrak{s}}}

\def\nn{\relax\ifmmode{\mathbb N_{0}}\else$\mathbb N_{0}$\fi}
\def\lra{\longrightarrow}

\begin{document}

\title[UNIFORM BEHAVIOUR OF FROBENIUS CLOSURES]{UNIFORM BEHAVIOUR
OF THE FROBENIUS CLOSURES OF IDEALS GENERATED BY REGULAR SEQUENCES}
\author{MORDECHAI KATZMAN}
\address{Department of Pure Mathematics,
University of Sheffield, Hicks Building, Sheffield S3 7RH, United Kingdom\\
{\it Fax number}: 0044-114-222-3769}
\email{M.Katzman@sheffield.ac.uk}
\author{RODNEY Y. SHARP}
\address{Department of Pure Mathematics,
University of Sheffield, Hicks Building, Sheffield S3 7RH, United Kingdom\\
{\it Fax number}: 0044-114-222-3769}
\email{R.Y.Sharp@sheffield.ac.uk}

\thanks{Both authors were partially supported by the
Engineering and Physical Sciences Research Council of the United
Kingdom.}


\date{\today}


\begin{abstract}
This paper is concerned with ideals in a commutative Noetherian
ring $R$ of prime characteristic. The main purpose is to show that
the Frobenius closures of certain ideals of $R$ generated by
regular sequences exhibit a desirable type of `uniform' behaviour.
The principal technical tool used is a result, proved by R.
Hartshorne and R. Speiser in the case where $R$ is local and
contains its residue field which is perfect, and subsequently
extended to all local rings of prime characteristic by G.
Lyubeznik, about a left module over the skew polynomial ring
$R[x,f]$ (associated to $R$ and the Frobenius homomorphism $f$, in
the indeterminate $x$) that is both $x$-torsion and Artinian over
$R$.
\end{abstract}

\maketitle

\setcounter{section}{-1}
\section{\bf Introduction}
\label{in}

Let $R$ be a commutative Noetherian ring of
prime characteristic $p$, and let $\fa$ be a proper ideal of $R$.
For $n \in \nn$ (we use $\nn$ (respectively $\N$) to denote the
set of non-negative (respectively positive) integers), the {\em
$n$-th Frobenius power\/} $\fa^{[p^n]}$ of $\fa$ is the ideal of
$R$ generated by all $p^n$-th powers of elements of $\fa$. Also
$R^{\circ}$ denotes the complement in $R$ of the union of the
minimal prime ideals of $R$.

An element $r \in R$ belongs to the {\em tight closure $\fa^*$ of
$\fa$\/} if and only if there exists $c \in R^{\circ}$ such that
$cr^{p^n} \in \fa^{[p^n]}$ for all $n \gg 0$. The theory of tight
closure was invented by M. Hochster and C. Huneke \cite{HocHun90},
and many applications have been found for the theory: see
\cite{Hunek96}.

The {\em Frobenius closure $\fa^F$ of $\fa$}, defined as
$$ \fa ^F := \big\{ r
\in R : \mbox{there exists~} n \in \nn \mbox{~such that~} r^{p^n}
\in \fa^{[p^n]}\big\},
$$
is another ideal relevant to the theory of tight closure. Since
$\fa^F$ is finitely generated, there exists $m_0 \in \nn$ such
that $(\fa^F)^{[p^{m_0}]} = \fa^{[p^{m_0}]}$, and we define
$Q(\fa)$ to be the smallest power of $p$ with this property. An
interesting question is whether the set $\{Q(\fb) : \fb \mbox{~is
a proper ideal of~} R\}$ of powers of $p$ is bounded. A simpler
question is whether, for a given proper ideal $\fa$ of $R$, the
set $\left\{Q\!\left(\fa^{[p^n]}\right) : n \in \nn\right\}$ is
bounded. The main result of this paper shows that the latter
question has an affirmative answer when $\fa$ is generated by a
regular sequence. The method of proof employs a result for local
$R$, about a left module over the skew polynomial ring $R[x,f]$
that is Artinian as an $R$-module, that was proved by R.
Hartshorne and R. Speiser \cite[Proposition 1.11]{HarSpe77} in the
case where the local ring $R$ contains its residue field which is
perfect, and subsequently extended to all local rings of
characteristic $p$ by G. Lyubeznik \cite[Proposition
4.4]{Lyube97}.  In Section \ref{cm} we apply this result, in the
case where $R$ is a Cohen--Macaulay local ring, to a top local
cohomology module viewed as a left $R[x,f]$-module. This is one
more instance where the Frobenius action on such a local
cohomology module, as described by K. E. Smith in \cite[3.2]{S94},
yields valuable insights (for other uses of this technique, see,
for example, \cite{F87}, \cite{FW87}, \cite{HW96},
\cite{HarSpe77}, \cite{Lyube97}).

The main result of the paper is presented in the final Section
\ref{rs}, where we use the modules of generalized fractions of the
second author and H. Zakeri \cite{32} to construct further modules
to which we can apply the Hartshorne--Speiser--Lyubezbik Theorem.

The uniform behaviour sought in this paper for Frobenius closures
has some similarity with the uniform behaviour of tight closures
that occurs when there exists a weak test element for $R$. A {\em
$p^{m_0}$-weak test element\/} for $R$ (where $m_0 \in \nn$) is an
element $c' \in R^{\circ}$ such that, for every ideal $\fb$ of $R$
and for $r \in R$, it is the case that $r \in \fb^*$ if and only
if $c'r^{p^n} \in \fb^{[p^n]}$ for all $n \geq m_0$. A $p^0$-weak
test element is called a {\em test element\/}. It is a result of
Hochster and Huneke \cite[Theorem (6.1)(b)]{HocHun94} that an
algebra of finite type over an excellent local ring of
characteristic $p$ has a $p^{m_0}$-weak test element for some $m_0
\in \nn$. To illustrate the relevance of the work in this paper to
the theory of tight closure, we establish the following lemma.

\begin{lem}
\label{in.1} Let $\fa$ be a proper ideal of the commutative
Noetherian ring $R$ of prime characteristic $p$. Assume that $m_0
\in \nn$ is such that either

\begin{enumerate}
\item there exists a $p^{m_0}$-weak test element $c'$ for $R$, or

\item $((\fa^{[p^n]})^F)^{[p^{m_0}]} =
(\fa^{[p^n]})^{[p^{m_0}]}$ for all $n \in \nn$.
\end{enumerate}
Then, for $r \in R$, we have $r \in \fa^*$ if and only if
there exists $c \in R^{\circ}$ such that
$cr^{p^n} \in (\fa^{[p^n]})^F$ for all $n \gg 0$.
\end{lem}

\begin{proof} Only one implication requires proof. Suppose that $r \in
R$, $n_0 \in \nn$, and $c \in R^{\circ}$ are such that $cr^{p^n}
\in (\fa^{[p^n]})^F$ for all $n \geq n_0$. We show that $r \in
\fa^*$.

Now $\fb^F \subseteq \fb^*$ for each ideal $\fb$ of $R$. Therefore, in case (i), we have
$$
c'(cr^{p^n})^{p^{m_0}} \in (\fa^{[p^n]})^{[p^{m_0}]} \quad \mbox{~for all~} n \geq n_0.
$$
In case (ii), we have
$$
(cr^{p^n})^{p^{m_0}} \in (\fa^{[p^n]})^{[p^{m_0}]} \quad \mbox{~for all~} n \geq n_0.
$$
If we let $\widetilde{c} = c'$ in case (i) and $\widetilde{c} = 1$
in case (ii), then we have (in both cases)
$(\widetilde{c}c^{p^{m_0}})r^{p^{n+m_0}} \in \fa^{[p^{n+m_0}]}$
for all $n \geq n_0$, and $\widetilde{c}c^{p^{m_0}} \in
R^{\circ}$. Hence $r \in \fa^*$.
\end{proof}

We also draw the reader's attention to the concept of {\em test exponent\/}
in tight closure theory introduced by Hochster and Huneke in
\cite[Definition 2.2]{HocHun00}. Let $c$ be a test element for a reduced
commutative Noetherian ring of prime characteristic $p$, and let $\fa$ be an
ideal of $R$. A test exponent for $c$, $\fa$ is a power $q = p^{e_0}$ (where
$e_0 \in \nn$) such that if, for an $r \in R$, we have $cr^{p^e} \in
\fa^{[p^e]}$ for {\em one single\/} $e \geq e_0$, then $r \in \fa^*$ (so that
$cr^{p^n} \in \fa^{[p^n]}$ for all $n \in \nn$). In \cite{HocHun00}, it is shown
that this concept has strong connections with the major open problem about whether
tight closure commutes with localization; indeed, to quote Hochster and Huneke,
`roughly speaking, $\ldots$ test exponents exist, in general, if and only if tight
closure commutes with localization'.

In \cite[Discussion 5.3]{HocHun00}, Hochster and Huneke raise the question as to
whether there might conceivably exist (when $R$ and $c$ satisfy certain
conditions) a `uniform test exponent' for $c$, that is, a power
of $p$ that is a test exponent for $c$, $\fb$ for {\em all\/} ideals $\fb$ of $R$
simultaneously. There are some similarities between this question and our question
(raised in the third paragraph of this Introduction) about whether the set
$\{Q(\fb) : \fb \mbox{~is
a proper ideal of~} R\}$ is bounded, and so we are hopeful that the work in this
paper might give some pointers for the major questions raised by Hochster and
Huneke.

\section{\bf Notation, terminology and the Hartshorne--Speiser--Lyubeznik Theorem}
\label{nt}

\begin{ntn}
\label{nt.1} Throughout the paper, $A$ will denote a general commutative Noetherian
ring and
$R$ will denote a commutative
Noetherian ring of prime characteristic $p$. For an ideal $\fc$ of $A$ and an
$A$-module $M$, we set $
\Gamma_{\fc}(M) := \left\{ m \in M : \mbox{there exists~} h \in
\N \mbox{~such that~} \fc^hm = 0 \right\}.
$

We shall always
denote by $f:R\lra R$ the Frobenius homomorphism, for which $f(r)
= r^p$ for all $r \in R$.  Throughout, $\fa$ will denote a general
proper ideal of $R$. We shall work with the
 skew polynomial ring $R[x,f]$ associated to $R$ and $f$
in the indeterminate $x$ over $R$. Recall that $R[x,f]$ is, as a
left $R$-module, freely generated by $(x^i)_{i \in \nn}$,
 and so consists
 of all polynomials $\sum_{i = 0}^n r_i x^i$, where  $n \in \nn$
 and  $r_0,\ldots,r_n \in R$; however, its multiplication is subject to the
 rule
 $$
  xr = f(r)x = r^px \quad \mbox{~for all~} r \in R\/.
 $$
\end{ntn}

\begin{lemdef}
\label{nt.2a} Let $Z$ be a left $R[x,f]$-module. Then the set
$$\Gamma_x(Z) := \left\{ z \in Z : x^jz = 0 \mbox{~for some~} j
\in \N \right\}$$ is an $R[x,f]$-submodule of $Z$, called the\/
{\em $x$-torsion submodule} of $Z$.
\end{lemdef}

\begin{proof} For $j \in \N$, $z \in Z$ and $r \in R$, we have
$x^jrz = r^{p^j}x^jz$; the claim follows easily.
\end{proof}

The following lemma enables one to see quickly that, in certain
circumstances, an $R$-module $M$ has a structure as left
$R[x,f]$-module extending its $R$-module structure.

\begin{lem}
\label{nt.3} Let $G$ be an $R$-module and let $\xi : G \lra G$ be a
$\Z$-endomorphism of $G$ such that $\xi(rg) = r^p\xi (g)$ for all $r \in R$ and
$g \in G$. Then the $R$-module structure on $G$ can be extended to a structure
of left $R[x,f]$-module in such a way that $xg = \xi (g)$ for all $g \in G$.
\end{lem}

\begin{proof} Note that, if we denote by $\mu_r$, for $r \in R$, the
$\Z$-endomorphism of $G$ given by multiplication by $r$,
then $\xi \circ \mu_r = \mu_{r^p} \circ \xi$ for all $r \in R$. We
use $\Endo_{\Z}(G)$ to denote the ring of $\Z$-endomorphisms of
the Abelian group $G$. In view of the universal property of
the skew polynomial ring $R[x,f]$, the above shows that there is a ring
homomorphism $\phi : R[x,f] \lra \Endo_{\Z}(G)$ for which
$\phi(x) = \xi$ and $\phi (r) = \mu_r$ for all $r \in R$. The
claim is now immediate from the fact that $G$ has a natural
structure as a left module over $\Endo_{\Z}(G)$.
\end{proof}

Crucial to the work in this paper is the following extension, due to G. Lyubeznik, of
a result of R. Hartshorne and
R. Speiser. It shows that, when $R$ is local, an
$x$-torsion left $R[x,f]$-module which is Artinian
(that is, `cofinite' in the terminology of Hartshorne and Speiser)
as an $R$-module exhibits a certain uniformity of behaviour.

\begin{thm} [G. Lyubeznik {\cite[Proposition 4.4]{Lyube97}}]
\label{hs.4}  {\rm (Compare Hartshorne--Speiser \cite[Proposition
1.11]{HarSpe77}.)} Suppose that $(R,\fm)$ is local, and let $G$ be
a left $R[x,f]$-module which is Artinian as an $R$-module. Then
there exists $e \in \N_0$ such that $x^e\Gamma_x(G) = 0$.
\end{thm}

Hartshorne and Speiser first proved this result in the
case where $R$ is local and contains its residue field which is perfect.
Lyubeznik applied his theory of $F$-modules to obtain the result
without restriction on the local ring $R$ of characteristic $p$.

\begin{defi}
\label{hslno} Suppose that $(R,\fm)$ is local, and let $G$ be
a left $R[x,f]$-module which is Artinian as an $R$-module.  By the
Hartshorne--Speiser--Lyubeznik Theorem \ref{hs.4},
there exists $e \in \nn$ such that $x^e\Gamma_x(G) = 0$: we call
the smallest such $e$ the {\em Hartshorne--Speiser--Lyubeznik number\/},
or {\em HSL-number\/} for short, of $G$.
\end{defi}

\section{\bf Parameter ideals in Cohen--Macaulay local rings}
\label{cm}

In this section we give an example of the use the
Hartshorne--Speiser--Lyubeznik Theorem to obtain a result, of the
type mentioned in the Introduction, about uniform behaviour of
Frobenius closures in a Cohen--Macaulay local ring $(R,\fm)$ of
positive dimension $d$. The argument is based on the well-known
$R[x,f]$-module structure that is carried by the top local
cohomology module $H^d_{\fm}(R)$. While this $R[x,f]$-module
structure has been used by several authors in the past, the fact
that this structure is independent of the choice of a system of
parameters for $R$ has not always been transparently clear from
the earlier accounts offered; as this independence is crucial for
our work, we shall make some comments about it.

\begin{rmd}
\label{lc.1} In this section, we shall sometimes use $R'$ to
denote $R$ regarded as an $R$-module by means of $f$, at points
where such care can be helpful.  Let $i \in \nn$.

\begin{enumerate}
\item With this notation, $f : R \lra R'$ becomes a homomorphism
of $R$-modules, and so, for the ideal $\fa$ of $R$, it induces an
$R$-homomorphism $ H^i_{\fa}(f) : H^i_{\fa}(R) \lra H^i_{\fa}(R').
$ \item The Independence Theorem for local cohomology (see
\cite[4.2.1]{LC}) applied to the ring homomorphism $f : R \lra R$
yields an $R$-isomorphism $ \nu^i_{R}: H^i_{\fa}(R')
\stackrel{\cong}{\lra} H^i_{\fa^{[p]}}(R)$, where
$H^i_{\fa^{[p]}}(R)$ is regarded as an $R$-module via $f$. Since
$\fa$ and $\fa^{[p]}$ have the same radical, $H^i_{\fa}$ and
$H^i_{\fa^{[p]}}$ are the same functor.

\item It is important for an understanding of
this paper to note that the $R$-isomorphism $\nu^i_{R}$ does not
depend on any choice of generators for $\fa$ or, for that matter,
for any ideal having the same radical as $\fa$. Indeed,
$\nu^i_{R}$ is a constituent isomorphism in a natural equivalence
of functors $\nu^i$ that forms part of an isomorphism of connected
sequences of functors that is uniquely determined by the identity
natural equivalence from $\Gamma_{\fa}$ to $\Gamma_{\fa^{[p]}}$:
see \cite[4.2.1]{LC} for details.

\item Composition yields a $\Z$-endomorphism
$
\xi := \nu^i_{R} \circ H^i_{\fa}(f) : H^i_{\fa}(R) \lra H^i_{\fa}(R)
$
which is such that $\xi(r\gamma) = r^p\xi(\gamma)$ for all $\gamma \in H^i_{\fa}(R)$
and $r \in R$. It therefore follows from Lemma \ref{nt.3} that $H^i_{\fa}(R)$
has a natural structure as left $R[x,f]$-module in which $x\gamma = \xi(\gamma)$
for all $\gamma \in H^i_{\fa}(R)$.

We emphasize that this $R[x,f]$-module structure on $H^i_{\fa}(R)$
does not depend on any choice of generators for an ideal having
the same radical as $\fa$.
\end{enumerate}
\end{rmd}

We can now define an invariant of a local ring $R$ of characteristic $p$.

\begin{defi}
\label{re.2d} Suppose that $(R,\fm)$ is a $d$-dimensional local
ring. By \ref{lc.1}, the top local cohomology module
$H^d_{\fm}(R)$, which is well known to be Artinian as an
$R$-module, has a natural structure as a left $R[x,f]$-module. We
define $\eta(R)$ to be the HSL-number (see \ref{hslno}) of
$H^d_{\fm}(R)$.
\end{defi}

To exploit the properties of this invariant, we are going to use a
concrete description of a typical element of $H^d_{\fm}(R)$ and
the way in which the indeterminate $x$ acts on such an element.

\begin{disc}
\label{cm.r1} If one chooses generators $b_1, \ldots, b_s$ for an
ideal $\fb$ having the same radical as $\fa$, then one can
represent the local cohomology module $H^i_{\fa}(R)$ quite
concretely, either as a direct limit of homology modules of Koszul
complexes, as in Grothendieck \cite[Theorem 2.3]{Groth67} or
\cite[\S 5.2]{LC}, or as a cohomology module of a \u{C}ech
complex, as in \cite[\S 5.1]{LC}. These representations can lead
to an explicit formula for the effect of $\xi$ (of \ref{lc.1}(iv))
on an element of $H^i_{\fa}(R)$. The following illustration for
the top local cohomology module $H^s_{\fa}(R)$ with respect to
$\fa$ is relevant to the work in this paper.

\begin{enumerate}
\item Let $M$ be an $R$-module. We are going to use the description
of the local cohomology module $H^s_{\fa}(M)$ as the $s$-th
cohomology module of the \u{C}ech complex of $M$ with respect to
$b_1, \ldots, b_s$. Thus $H^s_{\fa}(M)$ can be represented as the
residue class module of $M_{b_1 \ldots b_s}$ modulo the image,
under the \u{C}ech `differentiation' map, of
$\bigoplus_{i=1}^sM_{b_1 \ldots b_{i-1}b_{i+1}\ldots b_s}$. See
\cite[\S 5.1]{LC}. We use `$\left[\phantom{=} \right]$' to denote
natural images of elements of $M_{b_1\ldots b_s}$ in this residue
class module.

Denote the product $b_1 \ldots b_s$ by $b$. A typical element of
$H^s_{\fa}(M)$ can be represented as $ \left[m/b^n\right]$ for
some $m \in M$ and $n \in \nn$; moreover, for $m, m_1 \in M$ and
$n, n_1 \in \nn$, we have $ \left[m/b^n\right] =
\left[m_1/b^{n_1}\right] $ if and only if there exists $k \in \nn$
such that $k \geq \max\{n,n_1\}$ and $ b^{k-n}m - b^{k-n_1}m_1 \in
(b_1^k, \ldots, b_s^k)M. $ In particular, it should be noted that,
if $b_1, \ldots, b_s$ form an $M$-sequence, then $
\left[m/b^n\right] = 0$ if and only if $m \in (b_1^n, \ldots,
b_s^n)M$, by \cite[Theorem 3.2]{O'Car83}, for example. \item For
an element $v \in R$, there is an $R$-isomorphism $\omega : R'_{v}
\stackrel{\cong}{\lra} R_{v^p}$, where $R_{v^p}$ is regarded as an
$R$-module via $f$, for which $\omega(r/v^n) = r/v^{np}$ for all
$r \in R'$ and $n \in \nn$. It is straightforward to use such
isomorphisms (and the uniqueness aspect of \cite[Theorem
4.2.1]{LC}) to see that, with the notation of (i) above and
\ref{lc.1}(ii),
$$
\nu^i_R \left(\left[\frac{r}{(b_1\ldots b_s)^n}\right]\right) =
\left[\frac{r}{(b_1\ldots b_s)^{np}}\right] \quad \mbox{~for all~}
r \in R'.
$$
\item It follows that the left $R[x,f]$-module structure on $H^s_{\fa}(R)$ is such
that
$$
x\left[\frac{r}{(b_1\ldots b_s)^n}\right] =
\left[\frac{r^p}{(b_1\ldots b_s)^{np}}\right] \quad \mbox{~for
all~} r \in R \mbox{~and~} n \in \nn.
$$
In this section, we shall use the above formula in the special
case in which $(R, \fm)$ is local, $\fa = \fm$, and $b_1, \ldots,
b_s$ is a system of parameters for $R$.
\end{enumerate}
\end{disc}

\begin{thm}
\label{re.3} Suppose that $(R,\fm)$ {\rm (}as in {\rm \ref{nt.1})}
is a $d$-dimensional Cohen--Macaulay local ring, where $d > 0$.
The invariant $\eta(R)$ of $R$ is defined in\/ {\rm \ref{re.2d}}.

Let $e \in \nn$. Then the following statements are equivalent:
\begin{enumerate}
\item $e \geq \eta(R)$;
\item $(\fq^F) ^{[p^e]} = \fq ^{[p^e]}$ for each
ideal $\fq$ of $R$ that can be generated by a full system of
parameters for $R$;
\item for one system of parameters $a_1, \ldots, a_d$ for
$R$, we have
\[ ((a_1^{n_1},
\ldots,a_d^{n_d})^F) ^{[p^e]} = (a_1^{n_1}, \ldots,a_d^{n_d})
^{[p^e]} \quad \mbox{~for all~} n_1, \ldots, n_d \in \N\mbox{;} \]
\item for one ideal $\fs$ of $R$ that can be generated by a full system of
parameters for $R$, we have
\[
((\fs^{[p^n]})^F) ^{[p^e]} = (\fs^{[p^n]})^{[p^e]}
\quad \mbox{~for all~} n \in \nn. \]
\end{enumerate}
\end{thm}

\begin{proof}
Since $R$ is Cohen--Macaulay, every system of parameters for $R$
forms an $R$-sequence.

(i) $\Rightarrow$ (ii) Let $b_1, \ldots, b_d$ be a system of
parameters for $R$, and let $\fq = (b_1, \ldots, b_d)R$. Let $r
\in \fq^F$. Thus there exists $n \in \nn$ such that $r^{p^n} \in
\fq^{[p^n]} = (b_1^{p^n}, \ldots, b_d^{p^n})R$. Use $b_1, \ldots,
b_d$ in the notation of \ref{cm.r1}(i) for $H^d_{\fm}(R)$, and
write $b := b_1\ldots b_d$. We have $ x^n\left[r/b\right] =
\left[r^{p^n}/b^{p^n}\right] = 0 $ in $H^d_{\fm}(R)$. It therefore
follows from the definition of $\eta(R)$ as the HSL-number of
$H^d_{\fm}(R)$ that $ x^e\left[r/b\right] =
x^{e-\eta(R)}x^{\eta(R)}\left[r/b\right] = 0, $ so that $
\left[r^{p^e}/b^{p^e}\right] = 0 $ and $r^{p^e} \in \fq^{[p^e]}$
because $b_1, \ldots, b_d$ form an $R$-sequence. Therefore
$(\fq^F) ^{[p^e]} = \fq ^{[p^e]}$.

(ii) $\Rightarrow$ (iii) This is immediate from the fact that, if
$a_1, \ldots, a_d$ is a system of parameters for $R$, then so too is
$a_1^{n_1}, \ldots,a_d^{n_d}$ for all $n_1, \ldots, n_d \in \N$.

(iii) $\Rightarrow$ (iv) This is clear.

(iv) $\Rightarrow$ (i) Suppose that $\fs$ is generated by the
system of parameters $a_1, \ldots, a_d$ for $R$. Use $a_1, \ldots,
a_d$ in the notation of \ref{cm.r1}(i) for $H^d_{\fm}(R)$, and
write $a := a_1\ldots a_d$. Let $\zeta \in H^d_{\fm}(R)$. There
exist $r \in R$ and $n \in \nn$ such that $ \zeta = \left[
r/a^{p^n}\right] $; furthermore, $\zeta \in
\Gamma_x(H^d_{\fm}(R))$ if and only if there exists $j \in \nn$
such that $ \left[r^{p^j}/(a^{p^n})^{p^j}\right] = x^j\zeta = 0, $
that is (since $a_1, \ldots, a_d$ form an $R$-sequence), if and
only if there exists $j \in \nn$ such that $r^{p^j} \in
(((a_1,\ldots, a_d)R)^{[p^n]})^{[p^j]}$. Therefore, $\zeta \in
\Gamma_x(H^d_{\fm}(R))$ if and only if $r \in (\fs^{[p^n]})^F$.
Since $((\fs^{[p^n]})^F) ^{[p^e]} = (\fs^{[p^n]})^{[p^e]}$ for all
$n \in \nn$, it follows that $x^e (\Gamma_x(H^d_{\fm}(R))) = 0$,
so that $e \geq \eta(R)$.
\end{proof}

As a special case of Theorem \ref{re.3}, we recover the result of R. Fedder
(see \cite[Proposition 1.4]{F87}) that $\fq^F = \fq$
for each ideal $\fq$ of $R$ (as in \ref{re.3})
that can be generated by a full system of parameters (that is, $R$ is $F$-contracted
in the sense of \cite[p.\ 49]{F87}) if and only if $\eta(R) = 0$ (that is,
$R$ is $F$-injective in the sense of \cite[Definition 1.7]{FW87}). We
point out, however, that Fedder also proved that, in order for these
conditions to be satisfied, it is sufficient that $(a_1, \ldots, a_d)^F =
(a_1, \ldots, a_d)$ for {\em one single\/} system of parameters
$a_1, \ldots, a_d$ for $R$.

Our final result in this section shows that the set of ideals
$\fa$ in a Cohen--Macaulay local ring $(R,\fm)$ (as in \ref{nt.1})
of positive dimension for which $(\fa^F) ^{[p^{\eta(R)}]} =
\fa^{[p^{\eta(R)}]}$ is larger than the set of all ideals
generated by full systems of parameters, as it contains all
parameter ideals.  (We use the term `parameter ideal' in the sense
of Smith \cite[Definition 2.8]{S94}; however, a proper ideal in a
Cohen--Macaulay local ring is a parameter ideal if and only if it
can be generated by part of a system of parameters.)

\begin{thm}
\label{re.7} Suppose that $(R,\fm)$ {\rm (}as in {\rm \ref{nt.1})}
is a $d$-dimensional Cohen--Macaulay local ring, where $d > 0$.
The invariant $\eta(R)$ of $R$ is defined in\/ {\rm \ref{re.2d}}.

Let $\fa$ be an ideal of $R$ generated by part of a system of
parameters. Then $(\fa^F) ^{[p^{\eta(R)}]} = \fa^{[p^{\eta(R)}]}$.
\end{thm}

\begin{proof}
In view of Theorem \ref{re.3}, we can, and do, assume that
$\height \fa < d$. There exist a system of parameters $a_1,
\ldots, a_d$ for $R$ and an integer $i \in \{0, \ldots, d-1\}$
such that $\fa = (a_1,\ldots, a_i)R$. Let $r \in \fa^F$. Then, for
each $n \in \N$, we have $r \in ((a_1,\ldots, a_i,a_{i+1}^n,
\ldots, a_d^n)R)^F$, and, since $a_1,\ldots, a_i,a_{i+1}^n,
\ldots, a_d^n$ is a system of parameters for $R$, it follows from
Theorem \ref{re.3} that $r^{p^{\eta(R)}} \in ((a_1,\ldots,
a_i,a_{i+1}^n, \ldots, a_d^n)R)^{[p^{\eta(R)}]}$. Therefore
\begin{align*}
r^{p^{\eta(R)}} & \in \bigcap_{n \in \N}
\left(a_1^{p^{\eta(R)}},\ldots, a_i^{p^{\eta(R)}},
a_{i+1}^{np^{\eta(R)}}, \ldots, a_d^{np^{\eta(R)}}\right)\!R\\ &
\subseteq \bigcap_{n \in \N} \left(\fa^{[p^{\eta(R)}]} +
\fm^{np^{\eta(R)}} \right) = \fa^{[p^{\eta(R)}]}
\end{align*}
by Krull's Intersection Theorem.  The result follows.
\end{proof}

\section{\bf Preparatory results about modules of generalized fractions}
\label{gf}

For terminology and notation concerning
modules of generalized fractions, the reader is referred to Sharp--Zakeri \cite{32}.
Our first few results in this section
concern the general commutative Noetherian ring $A$. For an $A$-module $M$, we
say that a sequence $a_1, \ldots, a_t$ of elements of $A$ is a {\em poor
$M$-sequence\/} precisely when $\left((a_1, \ldots, a_i)M :_M a_{i+1}\right) =
(a_1, \ldots, a_i)M$ for all $i = 0, \ldots, t-1$. (Thus a poor $A$-sequence
is just a {\em possibly improper regular sequence\/} in the sense of
\cite[Discussion (7.3)]{HocHun90}.)

\begin{rmk}
\label{cm.5} A helpful tool for working with modules of generalized
fractions, particularly in the context of this paper, is the Exactness Theorem,
which provides an exactness criterion for complexes of modules of
generalized fractions.

Let ${\mathcal U} = (U_i)_{i \in
\N}$ be a chain of triangular subsets on $A$ in the sense of
\cite[page 420]{O'Car83}, and let $M$ be an
$A$-module. By \cite[3.3]{35}, the complex of modules of
generalized fractions $C({\mathcal U}, M)$ is exact if and only
if, for each $i \in \N$, each member of $U_i$ is a poor
$M$-sequence. (Although the present second author and H. Zakeri
first proved this result, a shorter proof, which applies also in
the case in which the underlying commutative ring is not
necessarily Noetherian, was later provided by O'Carroll
\cite[3.1]{O'Car83}.)

For each $i
\in \N$, we shall let $W_i = W_i(A)$ denote the set of all poor
$A$-sequences of length $i$. It is an easy consequence of
\cite[Example 3.10]{32} that $W_i$ is a triangular subset of
$A^i$. In fact, ${\mathcal W} = {\mathcal W}(A) := (W_i(A))_{i \in
\N}$ is a chain of triangular subsets on $A$ in the sense of
\cite[p.\ 420]{O'Car83}, and so the above-mentioned Exactness Theorem yields
that the complex $C({\mathcal W}, A)$ is exact.
\end{rmk}

The following consequence of the Exactness Theorem for generalized
fractions is very useful.

\begin{prop} [Sharp--Zakeri {\cite[Corollary 3.15]{35}}]
\label{gf.2a} Let $n \in \N$, let $M$ be an $A$-module,
and let $U$ be a triangular subset of $A^n$ that consists entirely of poor
$M$-sequences. Then, for $m \in M$ and $(u_1, \ldots, u_n) \in U$,
it is the case that \[\frac{m}{(u_1, \ldots, u_n)} = 0 \mbox{~in~}
U^{-n}M \quad \mbox{~if and only if~} \quad m \in
\sum_{i=1}^{n-1}u_iM.\]
\end{prop}

\begin{thm}
\label{gf.1} Suppose that $(A,\fm)$ is a local ring of depth $t > 0$.
Then $$\Hom_A(A/\fm, (W_t\times \{1\})^{-(t+1)}A)$$ is a finitely
generated $A$-module.
\end{thm}

\begin{proof} We consider the complex $C({\mathcal W}(A), A)$ of
modules of generalized fractions.  Note that this is exact. It is
also a consequence of the Exactness Theorem for modules of
generalized fractions that there are exact sequences
$$
0 \lra A \stackrel{d^{0}}{\lra} W_1^{-1}A \stackrel{\pi_1}{\lra}
(W_1 \times \{1\})^{-2}A \lra 0
$$
and (for each $i \in \N$ with $i > 1$)
$$
0 \lra (W_{i-1} \times \{1\})^{-i}A \stackrel{e^{i-1}}{\lra}
W_i^{-i}A \stackrel{\pi_i}{\lra} (W_i \times \{1\})^{-(i+1)}A \lra
0,
$$
where $d^0(a) = a/(1)$ for all $a \in A$, $ \pi_{1}(a/(r_1)) =
a/(r_1,1)$ for all $a \in A$, $(r_1) \in W_{1}$,
$$
e^{i-1}\left(\frac{a}{(r_1,\ldots,r_{i-1},1)}\right) =
\frac{a}{(r_1,\ldots,r_{i-1},1)} \quad \mbox{~for all~} a \in A,
(r_1,\ldots,r_{i-1}) \in W_{i-1},
$$
and
$$
\pi_{i}\left(\frac{a}{(r_1,\ldots,r_i)}\right) =
\frac{a}{(r_1,\ldots,r_i,1)} \quad \mbox{~for all~} a \in A,
(r_1,\ldots,r_{i}) \in W_{i}.
$$
(It follows from Proposition \ref{gf.2a} that $d^0$ and the $e^j~(j\in \nn)$ are
monomorphisms.)

We next show that, for each $i \in \N$ with $i \leq t$ and all $j \in \nn$, we
have $\Ext_A^j(A/\fm, W_i^{-i}A) = 0$. This is easy when
$i = 1$, and so we suppose that $i > 1$. For each ${\mathbf{r}} := (r_1,
\ldots, r_i) \in W_i$, let
$$
U_{\mathbf{r}} := \left\{ (r_1^{n_1}, \ldots, r_i^{n_i}): n_k \in \N
\mbox{~for all~} k = 1, \ldots, i \right\}.
$$
This is a triangular subset of $A^i$, and it follows from
\cite[page 39]{32} that $\left(U_{\mathbf{r}}^{-i}A\right)_{{\mathbf{r}} \in W_i}$ can
be turned into a direct system in such a way that
$$
\lim_{\stackrel {\scriptstyle \longrightarrow}{{\mathbf{r}} \in W_i}}
U_{\mathbf{r}}^{-i}A \cong W_i^{-i}A.
$$
Since $A/\fm$ is a finitely generated $A$-module, we have
$$
\Ext_A^j(A/\fm, W_i^{-i}A) \cong \lim_{\stackrel {\scriptstyle
\longrightarrow}{{\mathbf{r}} \in W_i}} \Ext_A^j(A/\fm,U_{\mathbf{r}}^{-i}A).
$$
We note next that $\Ext_A^j(A/\fm,U_{\mathbf{r}}^{-i}A) = 0$ for each
$A$-sequence ${\mathbf{r}} := (r_1, \ldots, r_i) \in W_i$, because
multiplication by $r_i$ on $U_{\mathbf{r}}^{-i}A$ provides an automorphism,
by \cite[Lemma 2.1]{34}, so that, since $r_i \in \fm$,
multiplication by $r_i$ on $\Ext_A^j(A/\fm,U_{\mathbf{r}}^{-i}A)$ provides an
endomorphism that is both zero and an automorphism.

Next, if ${\mathbf{r}}' := (r'_1, \ldots, r'_i) \in W_i$ is such that
$U_{{\mathbf{r}}'}^{-i}A \neq 0$, then $r'_1, \ldots, r'_{i-1} \in \fm$ (by
\ref{gf.2a}), so that $(r'_1, \ldots, r'_{i-1})$ is actually an
$A$-sequence. If $r'_i \not\in \fm$, then (since $i-1 < t$) there
exists $s_i \in \fm$ such that $(r'_1, \ldots, r'_{i-1},s_i)$ is
an $A$-sequence; of course, $s_i \in Rr'_i = R$, and so
$$
(r'_1, \ldots, r'_{i-1},r'_i) \leq (r'_1, \ldots, r'_{i-1},s_i)
$$
in the partial order used to form the direct limit in
\cite[Proposition 3.5]{32}.

These considerations show that $\Ext_A^j(A/\fm, W_i^{-i}A) = 0$,
for all $i \in \N$ with $i \leq t$ and all $j \in \nn$.

It now follows from the long exact sequences that result from
application of the functor $\Hom_A(A/\fm,\: {\scriptscriptstyle
\bullet} \:)$ to the short exact sequences displayed in the first
paragraph of this proof that
\begin{align*}
\Hom_A(A/\fm,(W_t\times \{1\})^{-(t+1)}A)  &\cong
\Ext_A^1(A/\fm,(W_{t-1}\times \{1\})^{-t}A) \cong \cdots
\\&\cong \Ext_A^{t-1}(A/\fm,(W_{1}\times \{1\})^{-2}A) \\&\cong
\Ext_A^t(A/\fm, A),
\end{align*}
and, since the last module is finitely generated, the result
is proved.
\end{proof}

\begin{cor}
\label{gf.2} Suppose that $(A,\fm)$ is a local ring of positive
depth $t$. Then $\Gamma_{\fm}\left((W_t\times
\{1\})^{-(t+1)}A\right)$ is an Artinian $A$-module.
\end{cor}

\begin{proof}
Note that $$\Hom_A\big(A/\fm, \Gamma_{\fm}\big((W_t\times
\{1\})^{-(t+1)}A\big)\big) \cong \Hom_A(A/\fm, (W_t\times
\{1\})^{-(t+1)}A)\mbox{;}$$ the latter module is finitely
generated, by Theorem \ref{gf.1}. Thus $G :=
\Gamma_{\fm}\left((W_t\times \{1\})^{-(t+1)}A\right)$ has finitely
generated socle; since each element of $G$ is annihilated by some
power of $\fm$, it follows that $G$ is Artinian (by
E. Matlis \cite[Proposition 3]{MatliE60} or L. Melkersson
\cite[Theorem 1.3]{MelkeL90}, for example).
\end{proof}

So far, the work in this section has concerned the general commutative
Noetherian ring $A$.  We now specialize to the commutative Noetherian ring $R$
of characteristic $p$.

\begin{lem}
\label{gf.3} Let $n \in \N$ and let $U$ be a triangular subset of
$R^n$. Then the module of generalized fractions $U^{-n}R$ has a
structure as left $R[x,f]$-module with
$$
x\left(\frac{r}{(u_1, \ldots,u_n)}\right) = \frac{r^p}{(u_1^p,
\ldots,u_n^p)} \quad \mbox{~for all~} r \in R \mbox{~and~} (u_1,
\ldots,u_n) \in U.
$$
\end{lem}

\begin{proof} Suppose that $\mathbf{u} = (u_1,
\ldots,u_n)$, $\mathbf{v} = (v_1, \ldots,v_n)$ and $\mathbf{w} =
(w_1, \ldots,w_n)$ are three elements of $U$ for which there exist
lower triangular matrices $\mathbf{H}, \mathbf{K}$ with entries in
$R$ such that $ \mathbf{H}\mathbf{u}^T = \mathbf{w}^T =
\mathbf{K}\mathbf{v}^T $ and that $r, s \in R$ are such that
$|\mathbf{H}|r - |\mathbf{K}|s \in (w_1,\ldots,w_{n-1})R$. (As in
\cite{32}, we use $|\mathbf{H}|$ to denote the determinant of
$\mathbf{H}$, {\it etcetera}.)  If $\mathbf{G}$ is a matrix with
entries from $R$, we shall denote by $\mathbf{G}^p$ the matrix
having the same size as $\mathbf{G}$ obtained by raising each
entry of $\mathbf{G}$ to the $p$-th power. Application of the
Frobenius homomorphism yields that
$$
\mathbf{H}^p(\mathbf{u}^p)^T = (\mathbf{w}^p)^T =
\mathbf{K}^p(\mathbf{v}^p)^T
$$
and $|\mathbf{H}^p|r^p - |\mathbf{K}^p|s^p \in
(w_1^p,\ldots,w_{n-1}^p)R$. It follows that $ r^p/(u_1^p,
\ldots,u_n^p) = s^p/(v_1^p, \ldots,v_n^p) $ in $U^{-n}R$, and that
there is a $\Z$-endomorphism $\xi : U^{-n}R \lra U^{-n}R$ which is
such that
$$
\xi \left(\frac{r}{(u_1, \ldots,u_n)}\right) = \frac{r^p}{(u_1^p,
\ldots,u_n^p)} \quad \mbox{~for all~} r \in R \mbox{~and~} (u_1,
\ldots,u_n) \in U.
$$
Note that $\xi (s\alpha) = s^p\xi(\alpha)$ for all $s \in R$ and
$\alpha \in U^{-n}R$. The claim therefore follows from Lemma
\ref{nt.3}.
\end{proof}

\begin{prop}
\label{gf.4} Suppose that $(R,\fm)$ is a local ring of positive
depth $t$. Then $$\Gamma_{\fm}\!\left((W_t(R)\times
\{1\})^{-(t+1)}R\right)$$ is an $R[x,f]$-submodule of
$(W_t(R)\times \{1\})^{-(t+1)}R$ which is Artinian as an
$R$-module.
\end{prop}

\begin{proof} If $\alpha \in (W_t(R)\times
\{1\})^{-(t+1)}R$ is annihilated by $\fm^h$ for a positive integer
$h$, then $(\fm^{h})^{[p]}x\alpha = 0$, and so
$\Gamma_{\fm}\left((W_t(R)\times \{1\})^{-(t+1)}R\right)$ is an
$R[x,f]$-submodule of $(W_t(R)\times \{1\})^{-(t+1)}R$. The result
therefore follows from Corollary \ref{gf.2}.
\end{proof}

\section{\bf Families of ideals generated by regular sequences}
\label{rs}

The purpose of this final section of the paper is to establish a
theorem which yields, as a corollary, the promised result that
$\left\{Q\!\left(\fa^{[p^n]}\right) : n \in \nn\right\}$ is
bounded when $\fa$ is generated by a regular sequence.  Our main
theorem concerns a family $(\fb_{\lambda})_{\lambda\in \Lambda}$
of ideals of $R$ that can be generated by $R$-sequences with the
property that ${\mathcal P} := \bigcup_{\lambda\in \Lambda} \ass
\fb_{\lambda}$ is a finite set (we use $\ass \fb_{\lambda}$ or
$\ass (\fb_{\lambda})$ to denote $\Ass (R/\fb_{\lambda})$). We
believe that Proposition \ref{rs.1}, which shows that there is a
good supply of examples of such families, is well known, but as we
have been unable to locate a precise reference for it, we give a
brief indication of proof.

\begin{prop}
\label{rs.1} Let ${\mathbf{r}} := (r_1, \ldots, r_t)$ and
${\mathbf{s}} := (s_1, \ldots, s_t)$ be poor $A$-sequences.

\begin{enumerate}
\item If ${\mathbf{s}}A \subseteq {\mathbf{r}}A$, then $\ass
({\mathbf{r}}A) \subseteq \ass ({\mathbf{s}}A)$.

\item We have $\ass \left(r_1^{n_1}, \ldots, r_t^{n_t}\right)\!A =
\ass ({\mathbf{r}}A)$ for all $n_1, \ldots, n_t \in \N$.

\item If $\sqrt{{\mathbf{r}}A} = \sqrt{{\mathbf{s}}A}$, then
$\ass ({\mathbf{r}}A) = \ass ({\mathbf{s}}A)$.
\end{enumerate}
\end{prop}

\begin{proof} (i) There is a $t \times t$ matrix ${\mathbf{M}}$
with entries in $R$ such that ${\mathbf{s}}^T =
{\mathbf{M}}{\mathbf{r}}^T$. By O'Carroll \cite[Theorem
3.7]{O'Car84}, the map $A/{\mathbf{r}}A \longrightarrow
A/{\mathbf{s}}A$ induced by multiplication by the determinant of
${\mathbf{M}}$ is an $A$-monomorphism. Therefore $\ass
({\mathbf{r}}A) \subseteq \ass ({\mathbf{s}}A)$.

(ii) It is enough to establish this result in the case where $A$
is local and ${\mathbf{r}}A$ is a proper ideal. In a local ring
$A$, every permutation of an $A$-sequence is again an
$A$-sequence; it is therefore sufficient for us to establish the
result in the case where $n_1 = \cdots = n_{t-1} = 1$. This can be
achieved easily by use of part (i) in conjunction with an
inductive argument based on the exact sequence
$$
0 \lra A/(r_1, \ldots, r_{t-1},r_t)A \lra A/(r_1, \ldots,
r_{t-1},r_t^{n_t+1})A \lra A/(r_1, \ldots, r_{t-1},r_t^{n_t})A
\lra 0,
$$
in which the second homomorphism is induced by multiplication by $r_t^{n_t}$
and the third map is the canonical epimorphism.

(iii)
There exists $n \in \N$ such that $(r_1^n, \ldots, r_t^n)A \subseteq
{\mathbf{s}}A$. By parts (i) and (ii), we have
$$
\ass({\mathbf{s}}A) \subseteq \ass(r_1^n, \ldots, r_t^n)A =
\ass({\mathbf{r}}A).
$$ Now reverse the roles of $\mathbf{r}$ and $\mathbf{s}$ to complete the
proof.
\end{proof}

We are now ready to present the main theorem of the paper.

\begin{thm}
\label{re.1} Let $\Lambda$ be an indexing set and let
$(\fb_{\lambda})_{\lambda\in \Lambda}$ be a family of
ideals of $R$ that can be generated by $R$-sequences with the
property that ${\mathcal P} := \bigcup_{\lambda\in \Lambda} \ass
\fb_{\lambda}$ is a finite set. Then there exists $e \in
\nn$ such that $$(\fb_{\lambda}^F)^{[p^e]} =
\fb_{\lambda}^{[p^e]} \quad \mbox{~for all~} \lambda\in
\Lambda.$$
\end{thm}

\begin{proof} We can assume that no $\fb_{\lambda}$ is $0$.
Let $\fp \in {\mathcal P}$, and let $t_{\fp} = \depth R_{\fp}$. By
Proposition \ref{gf.4},
$$
\Gamma_{\fp R_{\fp}}\big( \left( W_{t_{\fp}}(R_{\fp}) \times
\{1/1\}\right)^{-(t_{\fp}+1)}R_{\fp}\big)
$$
is an $R_{\fp}[x,f]$-submodule of $\left( W_{t_{\fp}} \times
\{1/1\}\right)^{-(t_{\fp}+1)}R_{\fp}$ which is Artinian as an
$R_{\fp}$-module. Let $e_{\fp}$ be its HSL-number (see \ref{hslno}). Let
$$\left\{ \height \fp : \fp \in {\mathcal P}\right\} = \{h_1,
\ldots, h_w\}, \quad \mbox{~where~} h_1 < h_2 < \cdots < h_w.
$$

Noting that ${\mathcal P} := \bigcup_{\lambda\in \Lambda} \ass
\fb_{\lambda}$ is a finite set, we
define (for each $i = 1, \ldots, w$)
$$e_i := \max \left\{e_{\fp}: \fp \in {\mathcal P}
\mbox{~and~} \height \fp = h_i\right\},$$ and we claim that $e =
\sum_{i=1}^w e_i$ has the desired property.
However, before we embark on the proof of this claim, we make the
following observations, in which $\fb, \fp$ denote ideals of $R$ with $\fp$
prime, and $i \in \nn$.

\begin{enumerate}
\item We have $(\fb R_{\fp})^F = (\fb^F) R_{\fp}$ and $(\fb
R_{\fp})^{[p^i]} = (\fb^{[p^i]}) R_{\fp}$. \item If
$(\fb^F)^{[p^i]} = \fb^{[p^i]}$, then $(\fb^F)^{[p^{i+j}]} =
(\fb^F)^{[p^{i+j}]}$ for all $j \in \N$. \item If $(\fb^F)^{[p^i]}
\neq \fb^{[p^i]}$, then the $R$-module
$(\fb^F)^{[p^i]}/\fb^{[p^i]}$ has an associated prime ideal $\fq$.
This $\fq$ will be such that $((\fb^F)^{[p^i]})R_{\fq} \neq
(\fb^{[p^i]})R_{\fq}$, that is (in view of (i) above) such that
$((\fb R_{\fq})^F)^{[p^i]} \neq (\fb R_{\fq})^{[p^i]}$.  Note that
such a $\fq$ has to be an associated prime of $\fb^{[p^i]}$;
consequently, if $\fb$ can be generated by an $R$-sequence, then
it follows from Proposition \ref{rs.1}(ii) that such a $\fq$ has
to be an associated prime of $\fb$.
\end{enumerate}

We can now use these observations (i)--(iii) to see that, in order
to establish our claim that $e = \sum_{i=1}^w e_i$ has the desired
property, it is enough for us to show that, for each $\fp \in
{\mathcal P}$, it is the case that
$$
((\fb_{\lambda}R_{\fp})^F)^{[p^{e_1 + \cdots + e_i}]}/
(\fb_{\lambda}R_{\fp})^{[p^{e_1 + \cdots + e_i}]} = 0 \quad \mbox{~for all~}
\lambda\in \Lambda, \mbox{~where~} \height \fp = h_i.
$$
We suppose that this is not the case, and we let $\fp$ be a minimal
counterexample. Set $\height \fp = h_j$; there must exist $\mu \in \Lambda$
such that
$$
((\fb_{\mu}R_{\fp})^F)^{[p^{e_1 + \cdots + e_j}]}/
(\fb_{\mu}R_{\fp})^{[p^{e_1 + \cdots + e_j}]} \neq 0.
$$
Set $e' := \sum_{\gamma = 1}^{j-1}e_{\gamma}$ (interpreted as $0$ if $j = 1$).
By choice of $\fp$, each of the $R_{\fp}$-modules
$$
((\fb_{\mu}R_{\fp})^F)^{[p^{e'}]}/
(\fb_{\mu}R_{\fp})^{[p^{e'}]} \quad \mbox{~and~} \quad
((\fb_{\mu}R_{\fp})^F)^{[p^{e'+e_j}]}/
(\fb_{\mu}R_{\fp})^{[p^{e'+e_j}]}
$$
has $\fp R_{\fp}$ as its only possible associated prime, because a
smaller associated prime would lead to a contradiction to the
minimality of $\fp$. Therefore, both of the $R_{\fp}$-modules in
this last display have finite length.

Let $r_1, \ldots, r_t$
be an $R$-sequence that generates $\fb_{\mu}$. Note that $\fp \in \ass \fb_{\mu}$,
by point (iii) above; it follows that $t = \depth R_{\fp} = t_{\fp}$ and
$\fb_{\mu}R_{\fp}$ is generated by the maximal $R_{\fp}$-sequence
$r_1/1, \ldots, r_t/1$.

There exists $\rho
\in (\fb_{\mu}R_{\fp})^F$ such that $\rho^{p^{e'+e_j}} \not\in
(\fb_{\mu}R_{\fp})^{[p^{e'+e_j}]}$. Consider
$$\alpha :=
\frac{\rho}{{\displaystyle \left(\frac{r_1}{1}, \ldots,
\frac{r_t}{1},\frac{1}{1}\right)}} \in
\left( W_{t_{\fp}}(R_{\fp}) \times \{1/1\}\right)^{-(t_{\fp}+1)}R_{\fp}.
$$
We have
$$
x^{e'}\alpha = x^{e'}
\frac{\rho}{{\displaystyle \left(\frac{r_1}{1}, \ldots,
\frac{r_t}{1},\frac{1}{1}\right)}} =
\frac{\rho^{p^{e'}}}{{\displaystyle \left(\frac{r_1^{p^{e'}}}{1}, \ldots,
\frac{r_t^{p^{e'}}}{1},\frac{1}{1}\right)}}
\in \Gamma_{\fp R_{\fp}}\big(\left( W_{t_{\fp}}(R_{\fp}) \times
\{1/1\}\right)^{-(t_{\fp}+1)}R_{\fp}\big),
$$
by Proposition \ref{gf.2a} (because $\rho^{p^{e'}} \in
((\fb_{\mu}R_{\fp})^F)^{[p^{e'}]}$ and
$((\fb_{\mu}R_{\fp})^F)^{[p^{e'}]}/ (\fb_{\mu}R_{\fp})^{[p^{e'}]}$
has finite length). Also, $\alpha \in \Gamma_{x}\big(\left(
W_{t_{\fp}}(R_{\fp}) \times
\{1/1\}\right)^{-(t_{\fp}+1)}R_{\fp}\big)$ because $\rho \in
(\fb_{\mu}R_{\fp})^F$. However, it follows from Proposition
\ref{gf.2a} that $x^{e_j}(x^{e'}\alpha) \neq 0$, because
$\rho^{p^{e'+e_j}} \not\in (\fb_{\mu}R_{\fp})^{[p^{e'+e_j}]}$.
Hence
$$
x^{e'}\alpha \in \Gamma_{x}\left(\Gamma_{\fp R_{\fp}}\big(\left(
W_{t_{\fp}}(R_{\fp}) \times
\{1/1\}\right)^{-(t_{\fp}+1)}R_{\fp}\big)\right) \quad
\mbox{~but~} \quad x^{e_j}(x^{e'}\alpha) \neq 0.
$$
This is a contradiction, and so the theorem is proved.
\end{proof}

Corollary \ref{re.2} below is a special case of Theorem \ref{re.1}; it presents one
of the results mentioned in the Introduction.

\begin{cor}
\label{re.2} Suppose that the ideal $\fa$ of $R$ can be generated
by an $R$-sequence $r_1, \ldots, r_t$. Then there exists $e \in
\nn$ such that $((r_1^{n_1}, \ldots, r_t^{n_t})^F)^{[p^e]}
= (r_1^{n_1}, \ldots, r_t^{n_t})^{[p^e]}$ for all $n_1, \ldots,
n_t \in \N$. In particular, $((\fa^{[p^n]})^F) ^{[p^e]}
= (\fa^{[p^n]})^{[p^e]}$ for all $n \in \nn$.
\end{cor}

\begin{proof} For
all $n_1, \ldots, n_t \in \N$, the sequence $r_1^{n_1}, \ldots,
r_t^{n_t}$ is an $R$-sequence, and $\ass (r_1^{n_1}, \ldots,
r_t^{n_t})R = \ass \fa$ by Proposition \ref{rs.1}(ii). The result is
therefore immediate from Theorem \ref{re.1}.
\end{proof}

\begin{cor}
\label{rs.4} Let $\Lambda$ be an indexing set and let
$(\fb_{\lambda})_{\lambda\in \Lambda}$ be a family of
ideals of $R$ that can be generated by $R$-sequences with the
property that
$
{\mathcal Q} := \left\{ \sqrt{\fb_{\lambda}} :
\lambda\in \Lambda \right\}
$
is a finite set. Then there exists $e \in
\nn$ such that $$(\fb_{\lambda}^F)^{[p^e]} =
\fb_{\lambda}^{[p^e]} \quad \mbox{~for all~} \lambda\in
\Lambda.$$
\end{cor}

\begin{proof}
For each $\fq \in {\mathcal Q}$, we can select a $\lambda(\fq) \in
\Lambda$ such that $\sqrt{\fb_{\lambda(\fq)}} = \fq$. Since
${\mathcal Q}$ is finite, the set ${\mathcal P}' := \bigcup_{\fq
\in {\mathcal Q}} \ass (\fb_{\lambda(\fq)})$ is finite. Let
$\mu\in \Lambda$. Then $\sqrt{\fb_{\mu}}$ is equal to some $\fq
\in {\mathcal Q}$, and so $\sqrt{\fb_{\mu}} =
\sqrt{\fb_{\lambda(\fq)}}$. By Proposition \ref{rs.1}(iii),
$\ass(\fb_{\mu}) = \ass(\fb_{\lambda(\fq)}) \subseteq {\mathcal
P}'$. Hence $\bigcup_{\mu\in \Lambda} \ass \fb_{\mu}$ is a subset
of the finite set ${\mathcal P}'$, and so the result follows from
Theorem \ref{re.1}.
\end{proof}

\bibliographystyle{amsplain}

\end{document}